\documentclass[12pt,twoside]{article}
\usepackage{amssymb,amsmath} 
\newtheorem{theorem}{Theorem}[section]

\newtheorem{lemma}[theorem]{Lemma}

\newtheorem{definition}[theorem]{Definition}

\numberwithin{equation}{section}
\newcommand \proof{\noindent{\bf Proof. \ }}
\newcommand \epf{\hfill $\square$ \par}

\title{Multiple solutions to a Caffarelli-Kohn-Nirenberg type equation with asymptotically linear term}

\begin{document}

\author{{Benjin Xuan}
\thanks{Supported by Grant 10101024 and 10371116 from the National Natural Science Foundation of China.}\\
{\it Department of Mathematics}\\
{\it University of Science and Technology of China}\\
{\it Universidad Nacional de Colombia}\\
{\it e-mail:wenyuanxbj@yahoo.com}}

\date{}
\maketitle

\begin{abstract}
In this paper, we study the existence of multiple solutions to a
Caffarelli-Kohn-Nirenberg type equation with asymptotically linear
term at infinity. In this case, the well-known
Ambrosetti-Rabinowtz type condition doesn't hold, hence it is
difficult to verify the classical (PS)$_c$ condition. To overcome
this difficulty, we use an equivalent version of Cerami's
condition, which allows the more general existence result.


\end{abstract}

\section{Introduction.}\label{intro}

In this paper, we shall investigate the existence of multiple solutions to the following  Caffarelli-Kohn-Nirenberg type equation with asymptotically linear term at infinity:
\begin{equation}
\label{eq1.1}
\left\{
\begin{array}{l}
-\mbox{div\,}(|x|^{-ap}|Du|^{p-2}Du)=|x|^{-(a+1)p+c}f(u),\mbox{ in } \Omega\\[2mm]
u=0,\  \mbox{on\,} \partial \Omega,
\end{array}
\right.
\end{equation}
where $\Omega\subset \mathbb{R}^n$ is an open bounded domain with $C^1$ boundary and $0\in \Omega$, $1<p<n,\ 0\leq a<\frac{n-p}p,\ c>0$ and $f(u)$ satisfies the following conditions:

(f$_1$) $f\in C(\mathbb{R},\ \mathbb{R}),\ f(0)=0,\ f(-t)=-f(t)$ for all $t\in \mathbb{R}$;

(f$_2$) $\lim\limits_{t\to 0}\dfrac{f(t)}{|t|^{p-2}t}=0$;

(f$_3$) (Asymptotically linear at infinity) $\lim\limits_{|t|\to \infty}\dfrac{f(t)}{|t|^{p-2}t}=l<\infty$;

(f$_3^\prime$) (Subcritical growth) Suppose $\lim\limits_{|t|\to \infty}\dfrac{f(t)}{t}=\infty$, and there exists some $r\in [p, p_*)$ with $p_*:=\dfrac{(n-(a+1)p+c)p}{n-(a+1)p}$ such that
$$
|f(t)|\leq C_1|t|^{p-1}+c_2|t|^{r-1} \mbox{ for some $C_1, C_2>0$};
$$

(f$_4$) $\dfrac{f(t)}{|t|^{p-2}t}$ is nondecreasing in $t\geq 0$;

(f$_5$) $\lim\limits_{|t|\to \infty}\{f(t)t-pF(t)\}=+\infty$.

For $a=0,\ c=p=2$, Brezis and Nirenberg \cite{BN} considered the
existence of positive solutions of problem (\ref{eq1.1}) in the
case where $f(u)$ is a lower-order perturbation of the critical
nonlinearity. They first showed that lower-order terms can reverse
the nonexistence and cause positive solutions to exist. After the
pioneering work \cite{BN}, there are many existence and
non-existence results of problem (\ref{eq1.1}) in many different
cases, e.g., Guedda and Veron \cite{GV} considered the quasilinear
elliptic equations involving critical Sobolev exponents; Zhu and
Yang \cite{ZY1, ZY2} also considered the quasilinear elliptic
equations involving critical Sobolev exponents in bounded or
unbounded domains. \cite{ABC} and \cite{HYX} considered the
combined effects of concave and convex nonlinearities in the
semilinear or quasilinear elliptic equations respectively.
\cite{GY} and \cite{XBJ} considered the existence of multiple
solutions to critical quasilinear equations with singular
symmetric coefficient or cylindrical symmetric coefficient. All of
the above results were based on, among some other structural
conditions, the following well-known Ambrosetti-Rabinowtz type
condition (cf. \cite{AR}), that is, there exists $\alpha>p$ such
that for $u\in \mathbb{R}$ and $|u|$ large enough, there holds
\begin{equation}
\label{eq1.01}
\alpha F(u)\leq uf(u),
\end{equation}
where $F(u)=\displaystyle\int_0^u f(t)\,dt$. The Ambrosetti-Rabinowtz type condition (\ref{eq1.01}) ensures the boundedness of the Palais-Smale sequence of the corresponding energy functional, then combining with other conditions, one can deduced the convergence of the Palais-Smale sequence. By simple calculation, it is easy to show that condition (\ref{eq1.01}) implies that
$$
F(u)\geq C |u|^\alpha,
$$
for some $C>0$ and $|u|$ large enough. Noting that $\alpha>p$, condition (\ref{eq1.01}) is not satisfied if $f(u)$ is asymptotically linear at infinity, i.e., $f(u)$ satisfies condition (f$_3$). This causes difficulty in proving the boundedness of the Palais-Smale sequence. Recently, Li and Zhou \cite{LZ} used an improved Mountain Pass Lemma with the usual Palais-Smale condition replaced by the Cerami weaker compactness condition (C) to obtain the existence of multiple solutions to problem (\ref{eq1.1}) in the non-singular case, i.e., $a=0,\ c=p$. In this paper, we shall extend the results in \cite{LZ} to the singular case, i.e., $0\leq a<\frac{n-p}p,\ c>0$.

\section{Preliminaries}
First, we recall an equivalent version of Cerami's condition as follows (cf. \cite{CM}):
\begin{definition}[Condition (C)]
\label{def1.1}
Let $X$ be a Banach space. $E\in C^1(X, \mathbb{R})$ is said to satisfy condition (C) at level $d\in \mathbb{R}$, if the following fact is true: any sequence $\{u_m\}\subset X$,  which satisfies
$$
E(u_m)\to d\ \mbox{and } (1+\|u_m\|)\|E^\prime(u_m)\|_{X^\prime} \to 0,\ \mbox{as } m\to \infty,
$$
possesses a convergent subsequence.
\end{definition}

As shown in \cite{BBF}, \cite{CM} and \cite{LZ}, under condition (C), there also hold the deformation lemma, Mountain Pass Lemma or saddle point theorem and Symmetric Mountain Pass Lemma. The following version of Symmetric Mountain Pass Lemma is due to Li and Zhou  \cite{LZ}:
Let $X$ be an infinite dimensional real Banach space, $E\in C^1(X, \mathbb{R})$, $\hat A_0=\{u\in X: \, E(u)\geq 0\}$,
$$
\Gamma^*=\{h:\, h(0)=0,\ h \mbox{ is an odd homeomorphism of $X$ and } h(B(0,1)) \subset \hat A_0\};
$$
$$
\begin{array}{ll}
\Gamma_m &=\{K\subset X:\, K \mbox{ is compact, symmetric with respect to the origin}, \\[2mm]
  & \ \  \ \ \ \ \ \mbox{ and for any $h\in \Gamma^*$, there holds } \gamma (K\cap h(\partial B(0,1))) \geq m \},
\end{array}
$$
where $\gamma(A)$ is the Krasnoselskii genus of $A$.
If $\Gamma_m\neq \emptyset$, define
$$
b_m=\inf_{K\subset \Gamma_m}\max_{u\in K}E(u).
$$
\begin{lemma}[Symmetric Mountain Pass Lemma]
\label{lem1.2} Let $e_1, e_2,\cdots, e_m, \cdots$ be linearly
independent in $X$, and $X_i=\mbox{span\,}\{e_1,\cdots e_i\},\
i=1,2, \cdots, m, \cdots$. Suppose that $E\in C^1(X, \mathbb{R})$
satisfies $E(0)=0$, $E(-u)=E(u)$, and condition (C) at all level
$d\geq 0$. Furthermore, there exists $\rho>0,\ \alpha>0$ such that
$E(u)>0$ in $B(0,\rho)\setminus \{0\}$ and $E|_{\partial B(0,
\rho)} \geq \alpha.$

Then, if $X_m\cap \hat A_0$ is bounded, then $\Gamma_m \neq \emptyset$ and $b_m\geq \alpha>0$ is critical. Moreover, if $X_{m+i}\cap \hat A_0$ is bounded for all $i=1, \cdots, l$, and
$$
b=b_{m+1}=\cdots =b_{m+l},
$$
then $\gamma(K_b)\geq l$, where $K_b=\{u\in X;\, E(u)=b,\ E^\prime(u)=0\}$.
if $X_m\cap \hat A_0$ is bounded for all $m$, then $E(u)$ possesses infinitely many critical values.
\end{lemma}

To apply Lemma \ref{lem1.2} to obtain the multiple results of problem (\ref{eq1.1}) with $a\geq 0$, we need the following weighted Sobolev-Hardy inequality due to Caffarelli, Kohn and Nirenberg \cite{CKN}, which is called the Caffarelli-Kohn-Nirenberg inequality. Let $1<p<n$. For all $u\in C_0^\infty(\mathbb{R}^n)$, there is a constant $C_{a,b}>0$ such that
\begin{equation}
\label{eq1.2}
\Big(\int_{\mathbb{R}^n}|x|^{-bq}|u|^{q}\,dx \Big)^{p/q}\leq C_{a,b}\int_{\mathbb{R}^n}|x|^{-ap}|Du|^{p}\,dx,
\end{equation}
where
\begin{equation}
\label{eq1.3}
-\infty< a<\frac{n-p}p,\ a\leq b\leq a+1,\ q=p^*(a,b)=\frac{np}{n-dp},\ d=1+a-b.
\end{equation}

Let $\Omega\subset \mathbb{R}^n$ is an open bounded domain with $C^1$ boundary and $0\in \Omega$, ${\cal D}_a^{1,p}(\Omega)$ be the completion of $C_0^\infty(\mathbb{R}^n)$, with respect to the norm $\|\cdot\|$ defined by
$$
\|u\|=\Big(\int_{\Omega}|x|^{-ap}|Du|^{p}\,dx \Big)^{1/p}.
$$
From the boundedness of $\Omega$ and the standard approximation
argument, it is easy to see that (\ref{eq1.2}) holds for any $u\in
{\cal D}_a^{1,p}(\Omega)$ in the sense:
\begin{equation}
\label{eq1.4}
\Big(\int_{\Omega}|x|^{-\alpha}|u|^{t}\,dx \Big)^{p/t}\leq C \int_{\Omega}|x|^{-ap}|Du|^{p}\,dx,
\end{equation}
for $1\leq t\leq \frac{np}{n-p},\ \alpha \leq (1+a)t+n(1-\frac tp)$, that is, the embedding ${\cal D}_a^{1,p}(\Omega) \hookrightarrow L^t(\Omega, |x|^{-\alpha})$ is continuous, where $L^t(\Omega, |x|^{-\alpha})$ is the weighted $L^r$ space with norm:
$$
\|u\|_{t, \alpha}:=\|u\|_{L^t(\Omega, |x|^{-\alpha})}=\Big( \int_{\Omega}|x|^{-\alpha}|u|^{t}\,dx\Big)^{1/t}.
$$

In fact, we have the following compact embedding result which is an extension of the classical Rellich-Kondrachov compactness theorem (cf. \cite{CC} for $p=2$ and  \cite{XB2} for the general case). For the convenience of readers, we include the proof here.
\begin{theorem}[Compact embedding theorem]
\label{thm1.1}
Suppose that $\Omega\subset \mathbb{R}^n$ is an open bounded domain with $C^1$ boundary and $0\in \Omega$, $1<p<n,\ -\infty< a<\frac{n-p}p$. The embedding ${\cal D}_a^{1,p}(\Omega) \hookrightarrow L^t(\Omega, |x|^{-\alpha})$ is compact if $1\leq t< \frac{np}{n-p},\ \alpha < (1+a)t+n(1-\frac tp)$.
\end{theorem}
\proof The continuity of the embedding is a direct consequence of
the Caffarelli-Kohn-Nirenberg inequality (\ref{eq1.2}) or
(\ref{eq1.4}). To prove the compactness, let $\{u_m\}$ be a
bounded sequence in ${\cal D}_a^{1,p}(\Omega)$. For any $\rho>0$
with $B_\rho(0)\subset \Omega$ is a ball centered at the origin
with radius $\rho$, there holds $\{u_m\}\subset
W^{1,p}(\Omega\setminus B_\rho(0))$. Then the classical
Rellich-Kondrachov compactness theorem guarantees the existence of
a convergent subsequence of $\{u_m\}$ in $L^t(\Omega\setminus
B_\rho(0))$. By taking a diagonal sequence, we can assume without
loss of generality that $\{u_m\}$ converges in
$L^t(\Omega\setminus B_\rho(0))$ for any $\rho>0$.

On the other hand, for any $1\leq t< \frac{np}{n-p}$, there exists a $b\in (a, a+1]$ such that $t<q=p^*(a, b)=\frac{np}{n-dp},\ d=1+a-b\in [0,\ 1)$. From the Caffarelli-Kohn-Nirenberg inequality (\ref{eq1.2}) or (\ref{eq1.4}), $\{u_m\}$ is also bounded in $L^q(\Omega, |x|^{-bq})$. By the H\"{o}der inequality, for any $\delta>0$, there holds
\begin{equation}
\label{eq2.001}
\begin{array}{ll}
\displaystyle \int_{|x|<\delta}|x|^{-\alpha}|u_m-u_j|^{t}\,dx & \leq \Big( \displaystyle \int_{|x|<\delta}|x|^{-(\alpha-bt)\frac q{q-t}}\,dx\Big)^{1-\frac tq}\\[3mm]
& \ \ \ \ \times
 \Big(\displaystyle \int_{\Omega}|x|^{-bt}|u_m-u_j|^{t}\,dx\Big)^{t/q}\\[3mm]
& \leq C \Big( \displaystyle \int_0^\delta r^{n-1-(\alpha-bt)\frac q{q-t}}\,dr\Big)^{1-\frac tq}\\[3mm]
& =C \delta ^{n-(\alpha-bt)\frac q{q-t}},
\end{array}
\end{equation}
where $C>0$ is a constant independent of $m$. Since $\alpha
< (1+a)t+n(1-\frac tp)$, there holds $n-(\alpha-bt)\frac
q{q-t}>0$. Therefore, for a given $\varepsilon>0$, we first fix
$\delta>0$ such that
$$
\int_{|x|<\delta}|x|^{-\alpha}|u_m-u_j|^{t}\,dx \leq \frac{\varepsilon}2, \ \forall\ m, j\in \mathbb{N}.
$$
Then we choose $N\in  \mathbb{N}$ such that
$$
\int_{\Omega\setminus B_\delta(0)}|x|^{-\alpha}|u_m-u_j|^{t}\,dx \leq C_\alpha \int_{\Omega\setminus B_\delta(0)} |u_m-u_j|^{t}\,dx \leq \frac{\varepsilon}2, \ \forall\ m, j\geq N,
$$
where $C_\alpha=\delta ^{-\alpha}$ if $\alpha\geq 0$ and $C_\alpha=(\mbox{diam\,}(\Omega) )^{-\alpha}$ if $\alpha< 0$. Thus
$$
\int_{\Omega}|x|^{-\alpha}|u_m-u_j|^{t}\,dx \leq \varepsilon , \ \forall\ m, j\geq N,
$$
that is, $\{u_m\}$ is a Cauchy sequence in $L^q(\Omega, |x|^{-bq})$.
\epf

On $X={\cal D}_a^{1,p}(\Omega)$, we define functional $E$ as
$$
E(u)=\frac1p\int_\Omega  |x|^{-ap} |D u|^p\,dx-\int_\Omega  |x|^{-(a+1)p+c} F(u)\,dx.
$$
By Theorem \ref{thm1.1} and assumptions (f$_1$)--(f$_3$) or (f$^\prime_3$), it is easy to show that $E$ is well-defined and of class $C^1(X, \mathbb{R})$, the weak solutions of problem (\ref{eq1.1}) is equivalent to the critical points of $E$. (f$_1$) implies that $0$ is a trivial solution to problem (\ref{eq1.1}).

In order to express our main theorem, we also need some results of the eigenvalue problem correspondent to problem (\ref{eq1.1}) in \cite{XB1}. Let us first recall the main results of \cite{XB1}. Consider the nonlinear eigenvalue problem:
\begin{equation}
\label{eq1.5}
\left\{
\begin{array}{l}
-\mbox{div\,}(|x|^{-ap}|Du|^{p-2}Du)=\lambda |x|^{-(a+1)p+c}|u|^{p-2}u,\mbox{ in } \Omega\\[2mm]
u= 0, \ \ \mbox{on } \partial\Omega,
\end{array}
\right.
\end{equation}
where $\Omega\subset \mathbb{R}^n$ is an open bounded domain with $C^1$ boundary and $0\in \Omega$, $1<p<n,\ 0\leq  a<\frac{n-p}p,\ c>0$.

Let us introduce the following functionals in ${\cal D}_a^{1,p}(\Omega)$:
$$
\Phi(u):= \int_{\Omega} |x|^{-ap}|Du|^{p}\,dx,\ \mbox{and } J(u):=\int_{\Omega} |x|^{-(a+1)p+c}|u|^p\,dx.
$$
For $c>0$, $J$ is well-defined. Furthermore, $\Phi, J\in C^1({\cal D}_a^{1,p}(\Omega),\mathbb{R})$, and a real value $\lambda$ is an eigenvalue of problem (\ref{eq1.5}) if and only if there exists $u\in {\cal D}_a^{1,p}(\Omega)\setminus\{0\}$ such that $\Phi^\prime(u)=\lambda J^\prime(u)$. At this point, let us introduce set
$$
{\cal M}:=\{u\in {\cal D}_a^{1,p}(\Omega)\ :\ J(u)=1 \}.
$$
Then ${\cal M}\neq \emptyset$ and ${\cal M}$ is a $C^1$ manifold in ${\cal D}_a^{1,p}(\Omega)$. It follows from the standard Lagrange multiples arguments that eigenvalues of (\ref{eq1.5}) correspond to critical values of $\Phi|_{{\cal M}}$. From Theorem \ref{thm1.1}, $\Phi$ satisfies the (PS) condition on ${\cal M}$. Thus a sequence of critical values of $\Phi|_{{\cal M}}$ comes from the Ljusternik-Schnirelman critical point theory on $C^1$ manifolds. Let $\gamma(A)$ denote the Krasnoselski's genus on ${\cal D}_a^{1,p}(\Omega)$ and for any $k\in \mathbb{N}$, set
$$
\Gamma_k:=\{A\subset {\cal M}\ :\ A \mbox{ is compact, symmetric and } \gamma(A)\geq k\}.
$$
Then values
\begin{equation}
\label{eq1.6}
\lambda_k:=\inf_{A\in\Gamma_k }\max_{u\in A} \Phi(u)
\end{equation}
are critical values and hence are eigenvalues of problem (\ref{eq1.5}). Moreover, $\lambda_1\leq \lambda_2\leq \cdots \leq   \lambda_k\leq \cdots \to+\infty$.

One can also define another sequence of critical values minimaxing $\Phi$ along a smaller family of symmetric subsets of ${\cal M}$. Let us denote by $S^k$ the unit sphere of $\mathbb{R}^{k+1}$ and
$$
{\cal O}(S^k, {\cal M}):= \{h\in C(S^k, {\cal M})\,:\, h \mbox{ is odd}\}.
$$
Then for any $k\in \mathbb{N}$, the value
\begin{equation}
\label{eq1.7}
\mu_k:= \inf_{h\in{\cal O}(S^{k-1}, {\cal M}) }\max_{t\in S^{k-1}} \Phi(h(t))
\end{equation}
is an eigenvalue of (\ref{eq1.5}). Moreover $\lambda_k\leq \mu_k$. This new sequence of eigenvalues was first introduced by \cite{DR} and later used in \cite{CM2} and \cite{CM} for $a=0, c=p$.


\section{Main theorem}
In this section, we shall prove the following main theorem:
\begin{theorem}
\label{thm3.1}
Assume that $f$ satisfies assumptions (f$_1$)-(f$_3$), and $\lambda_k$ is given by (\ref{eq1.6}). Then the following hold:
\begin{enumerate}
\item[(i)] If $l\in (\lambda_k, +\infty)$ is not an eigenvalue of problem (\ref{eq1.5}), then problem (\ref{eq1.1}) has at least $k$ pairs of nontrivial solutions in $X$;

\item[(ii)] Suppose that assumptions (f$_5$) is satisfied, then the conclusion of (i) holds even if $l$ is an eigenvalue of problem (\ref{eq1.5});

\item[(iii)] If $l=+\infty$, and assumptions (f$_3^\prime$)-(f$_4$) hold, then problem (\ref{eq1.1}) has infinitely many nontrivial solutions.
\end{enumerate}
\end{theorem}

We shall prove the above main theorem by verifying the assumptions in the Symmetric Mountain Pass Lemma-Lemma \ref{lem1.2}. First, we prove some properties of functional $E$.

\begin{lemma}
\label{lem3.1}
If assumptions (f$_1$) and (f$_4$) hold, then any sequence $\{u_m\}\subset X$ with $<E^\prime (u_m), u_m>\to 0$ as $m\to \infty$ possesses a subsequence, still denoted by $\{u_m\}$, such that
\begin{equation}
\label{eq3.1}
E(tu_m)\leq \frac{1+t^p}{mp}+E(u_m)
\end{equation}
holds for all $t\in \mathbb{R}, m\in \mathbb{N}$.
\end{lemma}
\noindent\textbf{Proof. 1. }From assumption $<E^\prime (u_m), u_m>\to 0$ as $m\to \infty$, up to a subsequence, for $m\in \mathbb{N}$, we may assume that
\begin{equation}
\label{eq3.2}
-\frac1m \leq <E^\prime (u_m), u_m>=\|u_m\|_X^p-\int_\Omega |x|^{-(a+1)p+c} f(u_m)u_m\,dx\leq \frac1m.
\end{equation}

\noindent\textbf{2. }For any $t\in \mathbb{R}, m\in \mathbb{N}$, there holds
\begin{equation}
\label{eq3.3}
E(tu_m)\leq \frac{t^p}{mp}+\int_\Omega |x|^{-(a+1)p+c} \{\frac1p f(u_m)u_m-F(u_m)\}\,dx.
\end{equation}

In fact, for any $t>0, m\in \mathbb{N}$, let
$$
h(t)=\frac1pt^p f(u_m)u_m-F(tu_m).
$$
Then, from assumptions (f$_1$) and (f$_4$), a simple calculation implies that
$$
\begin{array}{ll}
h^\prime(t)&=|t|^{p-2}tf(u_m)u_m-f(tu_m)u_m\\[2mm]
&=t^{p-2}t u_m^p\big(\dfrac{f(u_m)}{|u_m|^{p-2}u_m}-\dfrac{f(tu_m)}{t^{p-2}t|u_m|^{p-2} u_m} \big)\\[3mm]
&=\begin{cases} \geq 0, \mbox{ for }0<t\leq 1\\[2mm]
\leq 0, \mbox{ for } t\geq 1
\end{cases}
\end{array}
$$
and
\begin{equation}
\label{eq3.4}
h(t)\leq h(1), \mbox{ for all }t>0.
\end{equation}
Therefore, from (\ref{eq3.2}) and (\ref{eq3.4}), $t>0, m\in \mathbb{N}$, there holds
\begin{equation}
\label{eq3.5}
\begin{array}{ll}
E(tu_m)&=\dfrac1p t^p\|u_m\|_X^2-\displaystyle \int_\Omega F(tu_m)\,dx\\[3mm]
&\leq  \dfrac1p t^p\{\dfrac1m+\displaystyle  \int_\Omega  |x|^{-(a+1)p+c} f(u_m)u_m\,dx\}\\[3mm]
& \ \ \ \ \ \
-\displaystyle \int_\Omega |x|^{-(a+1)p+c} F(tu_m)\,dx\\[3mm]
&\leq \dfrac{t^p}{mp}+\displaystyle \int_\Omega  |x|^{-(a+1)p+c} \{\dfrac1p t^p f(u_m)u_m-F(tu_m)\}\,dx\\[3mm]
&\leq \dfrac{t^p}{mp}+\displaystyle \int_\Omega |x|^{-(a+1)p+c} \{\dfrac1p f(u_m)u_m-F(u_m)\}\,dx.
\end{array}
\end{equation}
Combining the oddness of $f$, (\ref{eq3.5}) implies that (\ref{eq3.3}) holds for any $t\in \mathbb{R}, m\in \mathbb{N}$.

\noindent\textbf{3. }For any $m\in \mathbb{N}$, there holds
\begin{equation}
\label{eq3.6}
\int_\Omega |x|^{-(a+1)p+c} \{\dfrac1p f(u_m)u_m-F(u_m)\}\,dx\leq \frac1{mp}+E(u_m).
\end{equation}

In fact, from (\ref{eq3.2}), there holds
$$
\begin{array}{ll}
E(u_m)&=\dfrac1p \|u_m\|_X^p-\displaystyle \int_\Omega |x|^{-(a+1)p+c} F(u_m)\,dx\\[3mm]
 &\geq \dfrac1p\{-\dfrac1m+\displaystyle  \int_\Omega |x|^{-(a+1)p+c} f(u_m)u_m\,dx \} -\displaystyle \int_\Omega |x|^{-(a+1)p+c} F(u_m)\,dx,
\end{array}
$$
that is, (\ref{eq3.6}) holds.

Combining (\ref{eq3.3}) and (\ref{eq3.6}) proves that (\ref{eq3.1}) holds.
\epf

\begin{lemma}
\label{lem3.2}
Let $d\geq 0$. Assume that $f$ satisfies assumptions (f$_1$) and (f$_2$). Then the following hold:
\begin{enumerate}
\item[i)] $E$ satisfies condition (C) at level $d\geq 0$ if assumption (f$_3$) holds and $l$ is not an eigenvalue of problem (\ref{eq1.5});

\item[ii)] $E$ satisfies condition (C) at level $d\geq 0$ if assumption (f$_3$) holds, $l$ is an eigenvalue of problem (\ref{eq1.5}) and assumption (f$_5$) holds;

\item[iii)] $E$ satisfies condition (C) at level $d\geq 0$ if assumptions (f$_3^\prime$) and (f$_4$) hold.
\end{enumerate}
\end{lemma}
\noindent\textbf{Proof.} Suppose that $\{u_m\}\subset X$ is a (C) sequence, that is, as $m\to \infty$, there hold
\begin{equation}
\label{eq3.7}
E(u_m)\to d\geq 0
\end{equation}
and
\begin{equation}
\label{eq3.8}
(1+\|u_m\|_X)\|E^\prime(u_m)\|_{X^\prime}\to 0.
\end{equation}
It is easy to see that (\ref{eq3.8}) implies that as $m\to \infty$, there hold
\begin{equation}
\label{eq3.9}
\|u_m\|_X^p-\int_\Omega |x|^{-(a+1)p+c} f(u_m)u_m\,dx=o(1)
\end{equation}
and
\begin{equation}
\label{eq3.10}
\int_\Omega |x|^{-ap} |D w_m|^{p-2} D w_m \cdot D v \,dx-
\int_\Omega |x|^{-(a+1)p+c} f(u_m)v\,dx=o(1),
\end{equation}
for all $v\in X$, where, and in what follows, $o(1)$ denotes any
quantity that tends to zero as $m\to \infty$.

From the compact embedding theorem \ref{thm1.1} and the fact that $f$ satisfies the subcritical growth condition, to show that $E$ satisfies condition (C) at level $d\geq 0$, it suffices to show the boundedness of (C) sequence $\{u_m\}$ in $X$ for each case.

\noindent\textbf{i). }Suppose that (f$_3$) holds and $l$ is not an eigenvalue of problem (\ref{eq1.5}). Suppose, by contradiction, that there exists a subsequence, still denoted by $\{u_m\}$, such that as $m\to \infty$, there holds
$$
\|u_m\|\to +\infty.
$$
Define
\begin{equation}
\label{eq3.010}
p_m(x)=\begin{cases} \dfrac{f(u_m(x))}{|u_m(x)|^{p-2} u_m(x)}, \ & \mbox{if } u_m(x)\neq 0\\[2mm]
    0, \ & \mbox{if } u_m(x)= 0.
\end{cases}
\end{equation}
Then from assumptions (f$_1$)--(f$_3$), there exists $M>0$ such that
\begin{equation}
\label{eq3.011}
0\leq p_m(x)\leq M, \ \forall \, x\in \Omega.
\end{equation}
Let
$$
w_m=\frac{u_m}{\|u_m\|}.
$$
Obviously, $\{w_m\}$ is bounded in $X$. Then from Lemma \ref{thm1.1}, without loss of generality, assume that there exists $w\in X$ such that as $m\to \infty$, there hold
\begin{equation}
\label{eq3.11}
w_m\rightharpoonup w \ \mbox{weakly in } X,
\end{equation}
\begin{equation}
\label{eq3.12}
w_m\to w \ \mbox{a.e. in } \Omega
\end{equation}
and
\begin{equation}
\label{eq3.13}
w_m\to w \ \mbox{strongly in } L^r(\Omega, |x|^{-(a+1)p+c}), \mbox{ if } 1\leq r< p_*.
\end{equation}

It is easy to show that $w\not\equiv 0$. In fact, if $w\equiv 0$, then from (\ref{eq3.9}), (\ref{eq3.011}), (\ref{eq3.13})
 and the definitions of $p_m$ and $w_m$, as $m\to \infty$, there holds
$$
\begin{array}{ll}
1& =\|w_m\|=\displaystyle \int_\Omega |x|^{-(a+1)p+c} p_m|w_m|^p\,dx+o(1)\\[3mm]
&\leq M \displaystyle \int_\Omega  |x|^{-(a+1)p+c}  |w_m|^p\,dx+o(1)\to 0,
\end{array}
$$
which is a contradiction.

From (\ref{eq3.011}), there exists $h\in L^\infty(\Omega)$ with $0\leq h\leq M$ such that, up to a subsequence, as $m\to \infty$, there holds
$$
p_m\rightharpoonup h \mbox{ weakly$^*$ in } L^\infty(\Omega).
$$
Then from (\ref{eq3.13}), there hold
\begin{equation}
\label{eq3.14}
p_m|w_m|^{p-2}w_m \rightharpoonup h|w|^{p-2}w \ \mbox{weakly in } L^{p^\prime}(\Omega,  |x|^{-(a+1)p+c} )
\end{equation}
where $p^\prime=p/(p-1)$, and
\begin{equation}
\label{eq3.15}
\int_\Omega  |x|^{-(a+1)p+c}  p_m|w_m|^p\,dx \to \int_\Omega  |x|^{-(a+1)p+c}  h|w|^p\,dx.
\end{equation}
On the other hand, from (\ref{eq3.9}) and (\ref{eq3.10}), there hold
\begin{equation}
\label{eq3.16}
\|w_m\|_X^p=\int_\Omega  |x|^{-(a+1)p+c}  p_m|w_m|^p\,dx +o(1)
\end{equation}
and
\begin{equation}
\label{eq3.17}
\int_\Omega  |x|^{-ap} |D w_m|^{p-2} D w_m \cdot D v\, dx-
\int_\Omega |x|^{-(a+1)p+c} p_m |w_m|^{p-2} w_m v\,dx=o(1),
\end{equation}
for all $v\in X$, Therefore, from (\ref{eq3.14})--(\ref{eq3.17}), there hold
\begin{equation}
\label{eq3.18}
\|w_m\|_X^p=\int_\Omega |x|^{-(a+1)p+c} h|w|^p\,dx +o(1)
\end{equation}
and
\begin{equation}
\label{eq3.19}
\int_\Omega |x|^{-ap} |D w_m|^{p-2} D w_m \cdot D v \,dx=
\int_\Omega |x|^{-(a+1)p+c} h |w|^{p-2} w v\,dx+o(1),
\end{equation}
for all $ v\in X$. It follows from (\ref{eq3.18}) and  (\ref{eq3.19}) that
$$
\int_\Omega |x|^{-ap}(|D w_m|^{p-2} D w_m-|D w|^{p-2} D w)\cdot (Dw_m-Dw)\,dx \to 0
$$
as $m\to \infty$, which gives that
$$
Dw_m \to  Dw \ \ \ \mbox{ in } L^p(\Omega,|x|^{-ap}).
$$
Thus $w$ satisfies
\begin{equation}
\label{eq3.20}
\int_\Omega |x|^{-ap} |D w|^{p-2} D w \cdot  Dv\,dx-
\int_\Omega |x|^{-(a+1)p+c}  h|w|^{p-2}w v\,dx=0,
\end{equation}
for all $ v\in X$. Let
$$
\begin{array}{ll}
&\Omega^+=\{x\in \Omega:\, w(x)>0\},\\[2mm]
&\Omega^0=\{x\in \Omega:\, w(x)=0\},\\[2mm]
&\Omega^-=\{x\in \Omega:\, w(x)<0\}.\\[2mm]
\end{array}
$$
Then $u_m(x)\to +\infty$ as $m\to \infty$ if $x\in \Omega^+$, and $u_m(x)\to -\infty$ as $m\to \infty$ if $x\in \Omega^-$. From assumption (f$_3$), $h(x)=l$ for all $x\in \Omega^+\cup \Omega^-$. Thus (\ref{eq3.20}) implies that $w$ satisfies
\begin{equation}
\label{eq3.21}
\int_\Omega |x|^{-ap} |D w|^{p-2} D w \cdot  Dv\,dx=
l\int_\Omega |x|^{-(a+1)p+c}  |w|^{p-2}w  v\,dx, \ \forall \, v\in X.
\end{equation}
This means that $l$ is an eigenvalue of problem (\ref{eq1.5}), which contradicts our assumption, so $\{u_m\}$ is bounded in $X$.

\noindent\textbf{ii). }Suppose $l$ is an eigenvalue of problem (\ref{eq1.5}), we need the additional assumption (f$_5$).

From assumption (f$_5$), there exists $T_0>0$ such that
$$
f(t)t-pF(t)\geq 0\ \mbox{ for all }|t|\geq T_0
$$
and there exists $C_0=C_0(T_0)>0$ such that
\begin{equation}
\label{eq3.22}
\int_{\{|u_m|\leq T_0\}}|x|^{-(a+1)p+c} [f(u_m)u_m-pF(u_m)]\,dx\geq -C_0.
\end{equation}
Furthermore, under assumptions ($f_1$)-($f_3$), there exists $M>0$ such that
\begin{equation}
\label{eq3.24}
|f(t)|\leq M|t|^{p-1},\ |F(t)|\leq \frac Mp|t|^p,\ \forall \, t\in R.
\end{equation}
Let $S>0$ be the best embedding constant such that
\begin{equation}
\label{eq3.25}
\big(\int_\Omega |x|^{-(a+1)p+c} |u|^{p_*}\,dx\big)^{p/p_*} \leq S \int_\Omega
|x|^{-ap} |D u|^p\,dx, \ \forall \,u\in X.
\end{equation}
Let $K=(dp+C_0)(2MS)^{(n-(a+1)p)/c+1}$ where $d$ is defined by (\ref{eq3.7}), $C_0$ by (\ref{eq3.22}), $M$ by (\ref{eq3.24}) and $S$ by (\ref{eq3.25}). From assumption ($f_5$), there exists a $T=T(K)>T_0>0$ such that
\begin{equation}
\label{eq3.26}
f(t)t-pF(t) \geq K,\ \mbox{for all } |t|>T.
\end{equation}
For the above $T>0$ and each $m\geq 1$, set
$$
A_m=\{(x,y)\in \Omega:\, |u_m(x,y)|\geq T\},\ B_m=\{(x,y)\in \Omega:\, |u_m(x,y)|\leq T\}.
$$
From estimates (\ref{eq3.22}), (\ref{eq3.7}), (\ref{eq3.9}) and (\ref{eq3.26}), there holds
\begin{equation}
\label{eq3.27}
\begin{array}{ll}
dp+o(1)&=\displaystyle \int_\Omega |x|^{-(a+1)p+c} [f(u_m)u_m-pF(u_m)]\,dx\\[3mm]
&\geq \displaystyle\int_{A_m}|x|^{-(a+1)p+c}  [f(u_m)u_m-2F(u_m)]\,dx-C_0\\[2mm]
&\geq KC(A_m)-C_0,
\end{array}
\end{equation}
where $C(A_m):= \displaystyle\int_{A_m} |x|^{-(a+1)p+c} \,dx$, noting that $a<(n-p)/p,\ c>0$.

On the other hand, for any fixed $r>p$, from (\ref{eq3.7}) and (\ref{eq3.9}), there holds
\begin{equation}
\label{eq3.28}
(\frac1p-\frac1r)\|u_m\|_X^p-\int_\Omega|x|^{-(a+1)p+c} [F(u_m)-\frac1r f(u_m)u_m]\,dx=d+o(1).
\end{equation}
Since $\Omega$ is bounded and $f(t)\in C(\mathbb{R},\mathbb{R})$, there exists a constant $C=C(\Omega, f, T)$ such that
\begin{equation}
\label{eq3.29}
|\int_{B_m} |x|^{-(a+1)p+c} [F(u_m)-\frac1r f(u_m)u_m]\,dx|\leq C, \ \mbox{for all }m\geq 1.
\end{equation}
Then, from (\ref{eq3.26})-(\ref{eq3.29}) and H\"{o}lder inequality, there holds
$$
\begin{array}{ll}
d+o(1)&\geq (\frac1p-\frac1r)\|u_m\|_X^p-C- \displaystyle\int_{A_m}|x|^{-(a+1)p+c} [F(u_m)-\frac1r f(u_m)u_m]\,dx\\[3mm]
&\geq (\frac1p-\frac1r)\|u_m\|_X^p-C-(\frac1p-\frac1r)\displaystyle\int_{A_m} |x|^{-(a+1)p+c}f(u_m)u_m\,dx\\[2mm]
&\geq (\frac1p-\frac1r)\|u_m\|_X^p-C-(\frac1p-\frac1r)M\displaystyle\int_{A_m} |x|^{-(a+1)p+c}|u_m|^p\,dx\\[2mm]
&\geq (\frac1p-\frac1r)\|u_m\|_X^p-C-(\frac 1p-\frac1r)M\|u_m\|_{p_*}^p(C(A_m))                                                                                                                                          ^{c/[n-(a+1)p+c]}\\[2mm]
&\geq (\frac1p-\frac1r)\|u_m\|_X^p-C-(\frac12-\frac1r)MS\|u_m\|_X^p\big(\frac{dp+C_0}K+o(1)\big)^{c/[n-(a+1)p+c}\\[2mm]
&\geq \frac12(\frac1p-\frac1r)\|u_m\|_X^p-C-(\frac1p-\frac1r)MS\|u_m\|_X^po(1),
\end{array}
$$
that is, $\{u_m\}$ is bounded in $X$.

\noindent\textbf{iii). }Suppose that assumptions ($f_3^\prime$) and ($f_4$) hold. To prove the boundedness of (C) sequence $\{u_m\}$, we need Lemma \ref{lem3.1}. Set
\begin{equation}
\label{eq3.30}
t_m=\frac{(2dp)^{1/p}}{\|u_m\|_X},\ w_m=t_mu_m=\frac{(2dp)^{1/p}\ u_m}{\|u_m\|_X}.
\end{equation}
Then $\|w_m\|_X=(2dp)^{1/p}$ and $\{w_m\}$ is bounded in $X$. Hence, up to a subsequence, we may assume that: there exists $w\in X$ such that estimates (\ref{eq3.11})-(\ref{eq3.13}) also hold in this case. If $\|u_m\|_X\to \infty$, we claim that
\begin{equation}
\label{eq3.32}
w(x)\not\equiv 0.
\end{equation}
In fact, if $w(x)\equiv 0$ in $\Omega$, then (\ref{eq3.30}) and estimates (\ref{eq3.11})-(\ref{eq3.13}) imply that
\begin{equation}
\label{eq3.33}
\int_\Omega |x|^{-(a+1)p+c}  F(w_m)\,dx\to 0 \mbox{ and  } E(w_m)=2dp+o(1).
\end{equation}
However, applying Lemma \ref{lem3.1} with $t=(2dp)^{1/p}/\|u_m\|_X
\to 0$, there holds
$$
E(w_m)\leq \frac{1+|t|^p}{mp}+E(u_m)\to d,
$$
which contradicts (\ref{eq3.33}), thus (\ref{eq3.32}) holds.

On the other hand, similar to case i), (\ref{eq3.16}) holds. Let
$$
\bar \Omega=\Omega\setminus \{x\in \Omega:\, w(x)=0\}.
$$
Then $|\bar \Omega|>0$ by (\ref{eq3.32}). From assumptions ($f_2$), ($f_3^\prime$) and ($f_4$), $p_m(x)\geq 0$ and $p_m(x)\to \infty$ as $m\to \infty$ in $\bar \Omega$, where $p_m$ is defined by (\ref{eq3.010}). Hence, from (\ref{eq3.16}), there holds
$$
\begin{array}{ll}
2dp&=\liminf\limits_{m\to \infty}\|w_m\|_X^p
=\liminf\limits_{m\to \infty}\displaystyle\int_\Omega |x|^{-(a+1)p+c}  p_m(x)|w_m|^p\,dx\\[3mm]
&\geq \liminf\limits_{m\to \infty}\displaystyle\int_{\bar \Omega}|x|^{-(a+1)p+c}  p_m(x)|w_m|^p\,dx\\[3mm]
&\geq \displaystyle\int_{\bar \Omega}\liminf\limits_{m\to \infty}|x|^{-(a+1)p+c}  p_m(x)|w_m|^p\,dx=\infty,
\end{array}
$$
which is a contradiction, thus $\|u_m\|\not\to \infty$, that is, up to a subsequence, $\{u_m\}$ is bounded in $X$.
\epf

\noindent\textbf{Proof of Theorem \ref{thm3.1}.} We shall prove this theorem by verifying the assumptions of the Symmetric Mountain Pass Lemma-Lemma \ref{lem1.2}.

\noindent\textbf{1.} There exists $\rho>0,\ \alpha>0$ such that $E(u)>0$ in $B(0,\rho)\setminus \{0\}$ and $E|_{\partial B(0, \rho)} \geq \alpha.$

In fact, in each case, assumptions (f$_1$)-(f$_3$) or (f$_3^\prime$) imply that $E\in C^1(X, \mathbb{R})$ and for any $\varepsilon>0$, there exists $C_{\varepsilon}>0$ such that, for all $t\in \mathbb{R}$, there holds
\begin{equation}
\label{eq3.34}
|f(t)|\leq \varepsilon |t|^{p-1}+C_{\varepsilon}|t|^{r-1},\ \ |F(t)|\leq \varepsilon |t|^p+C_{\varepsilon}|t|^{r},
\end{equation}
from which, it is easy to see that there exists $\rho>0,\ \alpha>0$ such that $E(u)>0$ in $B(0,\rho)\setminus \{0\}$ and $E|_{\partial B(0, \rho)} \geq \alpha$.

\noindent\textbf{2.} By the Symmetric Mountain Pass Lemma-Lemma
\ref{lem1.2}, to prove Theorem \ref{thm3.1}, it suffices to prove
that for any $k\geq 1$, there exists a $k$-dimensional subspace
$X_k$ of $X$ and $R_k>0$ such that
\begin{equation}
\label{eq3.35}
E(u)\leq 0,\ \mbox{for all } u\in X_k\setminus B_{R_k},
\end{equation}
where $B_{R_k}$ is the ball in $X$ with radius $R_k$.

First, we prove (\ref{eq3.35}) in the case $l\in (\lambda_k, +\infty)$. Since $l>\lambda_k$, there is $\varepsilon>0$ such that $l-\varepsilon>\lambda_k$. By the definition of $\lambda_k$, there exists a $k$-dimensional subspace $X_k$ of $X$ such that, for the above $\varepsilon>0$, there holds
\begin{equation}
\label{eq3.36}
\sup_{u\in X_k\setminus \{0\}}\dfrac{\Psi(u)}{J(u)}\leq \lambda_k+\frac\varepsilon 2 <l-\frac\varepsilon 2,
\end{equation}
that is,
\begin{equation}
\label{eq3.37}
\inf_{u\in X_k\setminus \{0\}}\dfrac{J(u)}{\Psi(u)}> \frac1{l-\varepsilon/ 2}.
\end{equation}

On the other hand, by assumption (f$_3$), there holds
\begin{equation}
\label{eq3.38}
\lim_{|u|\to \infty} \frac{F(u)}{|u|^p}=\frac lp.
\end{equation}
Then, for the above $\varepsilon>0$, there exists $M>0$ large enough such that
\begin{equation}
\label{eq3.39}
\frac{F(u)}{|u|^p}> \frac 1p (l-\frac \varepsilon 4),\ \ \mbox{for all } |u|>M.
\end{equation}
Therefore, if $u\in X_k$ with $\|u\|=R$, there holds
\begin{equation}
\label{eq3.40}
\begin{array}{ll}
E(u)&= \dfrac1p R^p-\displaystyle \int_{\Omega} |x|^{-(a+1)p+c} F(u)\,dx\\[3mm]
&\leq  \dfrac1p R^2-\displaystyle \int_{|u|>M} |x|^{-(a+1)p+c} F(u)\,dx-C(M, \Omega)\\[3mm]
& \leq  \dfrac1p R^2-\dfrac 12 (l-\frac \varepsilon 4)\displaystyle \int_{\Omega}|x|^{-(a+1)p+c} |u|^p\,dx-C(M, \Omega)\\[3mm]
&= \dfrac {R^2}p \big[1- (l-\dfrac \varepsilon 4)\displaystyle \int_{\Omega}|x|^{-(a+1)p+c} \dfrac{|u|^p}{R^p}\,dx \big] -C(M, \Omega)\\[3mm]
&\leq  \dfrac {R^p}p \big(1- \dfrac{l-\varepsilon /4}{l- \varepsilon/ 2} \big)-C(M, \Omega)\\[3mm]
&<0,
\end{array}
\end{equation}
if $R\geq R_k$ and $R_k>0$ large enough.

If $l=+\infty$, similarly to (\ref{eq3.37}), for any $k\geq 1$, there exists $X_k$ such that
$$
\inf_{u\in X_k\setminus \{0\}}\dfrac{J(u)}{\Psi(u)} > \frac1{\lambda_k +1/ 2},
$$
and, similarly to (\ref{eq3.39}), it follows from assumption (f$_3^\prime$) and $l=+\infty$ that there exists $M_k>0$ such that
$$
\frac{F(u)}{|u|^p}> \frac 12 (\lambda_k +1),\ \ \mbox{for all } |u|>M_k.
$$
Then, if $u\in X_k$ with $\|u\|=R$, there holds
$$
E(u)\leq  \dfrac {R^p}p\big(1- \dfrac{\lambda_k +1}{\lambda_k +1/ 2} \big)-C(M_k, k,  \Omega)<0,
$$
if $R\geq R_k$ and $R_k>0$ large enough.
\epf


\noindent\textsc{Benjin Xuan} \\
Department of Mathematics\\
Universidad Nacional de Colombia\\
Bogot\'a, Colombia\\
and \\
University of Science and Technology of China\\
Hefei, Anhui, 230026\\
e-mail: bjxuan@matematicas.unal.edu.co, wenyuanxbj@yahoo.com

\end{document}